\documentclass[12pt]{amsart}
\usepackage{amsmath}
\usepackage{amsthm}
\usepackage{latexsym}
\usepackage{amssymb}
\newtheorem{theorem}{Theorem}
\newtheorem{proposition}{Proposition}

\newtheorem{corollary}{Corollary}
\newtheorem{definition}{Definition}

\newcounter{obsctr}

\newtheorem{remark}{Remark}

\renewcommand{\theequation}{\thesection.\arabic{equation}}
\newcount\madechanges
\def\caution#1{\ifnum \madechanges=1 \affixmessage{#1}%
\else \relax \fi}
\def\affixmessage#1{\marginpar{{\footnotesize \em #1} \openup
    -.3\baselineskip }}  
\begin{document}
\baselineskip 16pt
\def\A {{\mathcal{A}}}
\def\D {{\mathcal{D}}}
\def\R {{\mathbb{R}}}
\def\N {{\mathbb{N}}}
\def\C {{\mathbb{C}}}
\def\Z {{\mathbb{Z}}}
\def\l {\ell}
\def\ml {multline}
\def\multiline {\multline}
\def\lessim {\lesssim}
%
%
%
% the following are for better readability
%
%
%
\def\phi{\varphi}
\def\epsilon{\varepsilon}
\title{Global Analytic Hypoellipticity for a Class of
Quasilinear Sums of Squares of Vector Fields}    
\author{Makhlouf Derridj}
\address{UPRESA, UniversitŽ de Haute-Normandie, 
Laboratoire Raphael Salem; 
UMR 6085, CNRS, Cite Colbert, 76821
Mont-Saint-Aignan, FRANCE}
\email{Makhlouf.Derridj@univ-rouen.fr}
\author{David S. Tartakoff}
\address{Department of Mathematics, University
of Illinois at Chicago, m/c 247, 851 S.
Morgan St., Chicago IL  60607}
\email{dst@uic.edu}
%\date{\today}
\begin{abstract}
We prove a global analytic regularity result for an operator of
H\"ormander rank $2$ constructed from $2n$ quasilinear vector
fields on a compact manifold in $\R^{2n+1}.$
\end{abstract}
\maketitle
\pagestyle{myheadings}
\markboth{M. Derridj and D.S.Tartakoff}{Globally Analytic
Nonlinear Sums of Squares}
%
%%%%%%%%%%%%%%%%%%%%%%%%%%%%%%%%%%%%%%%%%%%
%%%%%%%%%%%%%%%%%%%%%%%%%%%
%%%%%%%%%%%%%%%%%%%% required by fancyheadings package 
%%%%%%%%%%%%%%%
%%%%%%%%%%%%%%%%%%%%%%%%%%%%%%%%%%%%%%%%%%%
%%%%%%%%%%%%%%%%%%%%%%%%%%%
%\thispagestyle{empty}
%\fancyhead{}
%\fancyhead[RE]{\small M. Derridj and D.S. Tartakoff}
%\fancyhead[LE]{\thepage}
%\fancyhead[LO]{\small Globally Analytic Nonlinear Sums of
%Squares}
%\fancyhead[RO]{\thepage}
%\fancyfoot[CO,CE]{}
%\renewcommand{\headrulewidth}{0pt}
%\renewcommand{\footrulewidth}{0pt}
%%%%%%%%%%%%%%%%%%%%%%%%%%%%%%%%%%%%%%%%%%%
%%%%%%%%%%%%%%%%%%%%%%%%%%%
%%%%%%%%%%%%%%%%%%%%%%%%%%%%%%%%%%%%%%%%%%%
%%%%%%%%%%%%%%%%%%%%%%%%%%%
%%%%%%%%%%%%%%%%%%%%%%%%%%%%%%%%%%%%%%%%%%%
%%%%%%%%%%%%%%%%%%%%%%%%%%%
%
%
%
\section{Introduction}
\renewcommand{\theequation}{\thesection.\arabic{equation}}
\setcounter{equation}{0}
\setcounter{theorem}{0}
\setcounter{proposition}{0}  
\setcounter{lemma}{0}
\setcounter{corollary}{0} 
\setcounter{definition}{0}

Our aim, in this work, is to prove a global analytic
regularity result on a compact manifold $M$ for some
quasilinear equations. 

A model of such results is the following: in $C^2$
take a bounded domain $\Omega$ with strictly pseudo-convex real
analytic boundary, M. Then one has two globally defined
independent real, real analytic vector fields $X_1$ and
$X_2,$ namely the real and imaginary parts of a (globally 
defined) holomorphic vector field $L=X_1-iX_2$ tangent to $M.$ 

Take a function $u$ in $C^\infty(M)$ and consider
an analytic matrix function $H(x,t)$ defined on a
neighborhood of $\{(x,u(x)):x\in M\}$ in $M\times \C$
and set 
$$Y=H(X) \hbox{ i.e., } Y(x) = H(x,u(x))X(x)$$
so that one obtains two $C^\infty$ vector fields,
$Y_1$ and $Y_2$ on $M.$

We assume that $H(x,t)$ is invertible so
that one can express $X=H^{-1}(Y)$ with $H^{-1}\in
C^\omega.$

Consider the operator 
$$P_u = Y_1^2+Y_2^2+a_1Y_1+a_2Y_2+b$$
with $a_j,b$ analytic and assume that $P_uu \in \A(M).$
Can one conclude that $u$ is analytic on $M$ if the
associated Levi form is non-degenerate? Note that
$P_uu=f$ is a quasi-linear equation.

The question is global. There are known local results
for more special equations (cf. \cite{Tartakoff-Zanghirati2003}) 

In higher dimensions, one generally does not have
globally defined vector fields $X_1, \ldots, X_{2n}$
related to a CR structure on $M$ induced by the complex
structure on $\C^n.$ But one can consider a (finite)
family of open sets $\{V_\ell\}_{1\leq\ell\leq p}$
covering $M$ and analytic vector fields
$\{X_{k,\ell}\}_{k=1,\ldots,2n}$ on $V_\ell$. Then we
may consider
$$Y_{(\ell)} = H(X_{(\ell)}) \hbox{ where } X_{(\ell)}
= \left(\begin{matrix}X_{1,\ell}\\ \vdots \\
X_{2n,\ell}\end{matrix}\right)$$
and the associated operator
$$P_{\ell,u} = \sum_{j=1}^{2n}Y_{j,\ell}^2 +
a_{j,\ell}Y_{j,\ell} + b_\ell, \hbox{ with }
a_{j,\ell}, b_\ell \hbox{ analytic }.$$

Now assume that for all $\ell,$ 
$$P_{\ell,u}u \in \A(V_\ell).$$

Then the question is: is $u$ analytic on $M$ under a
non-degeneracy hypothesis on the associated Levi
form?
\begin{theorem}
Under the above hypotheses, if $Pu$ is real analytic globally on
$M$ then so is any (moderately smooth) solution $u.$
\end{theorem}

A more interesting problem (as related to the
complex Laplacian on forms) is to consider a system.
But in this paper we consider only the scalar case.
Note that from results on $C^\infty$ regularity (cf. Xu
(), et al.), one need only assume that $u$ is in $C^{2,\alpha}$ in
our theorem. 

\section{Some notation and definitions}
\renewcommand{\theequation}{\thesection.\arabic{equation}}
\setcounter{equation}{0}
\setcounter{theorem}{0}
\setcounter{proposition}{0}  
\setcounter{lemma}{0}
\setcounter{corollary}{0} 
\setcounter{definition}{0}

Let $M$ be a compact, real analytic manifold of
dimension $2n+1\geq 3,$ and let $(V_j)_{j=1\ldots
p}$ be a covering of $M$ such that, in each $V_j,$
there are given $2n$ real analytic, real vector
fields $X_{1,j}, \ldots , X_{2n,j}$ such that 

\begin{itemize}
\item On $V_j\cap V_k$ every $X_{\ell,j}$ (resp.
$X_{\ell,k}$) is a linear combination of the
$(X_{\ell,k})_{\ell=1\ldots p}$ (resp. of the
$(X_{\ell,j})_{\ell=1\ldots p}$) with real analytic
coefficients. 
\item There exists a globally defined, real analytic real vector
field $T$ such that
$(X_{1,j}.
\ldots , X_{2n,j}, T)$ is a basis in $V_j$ and if 
\begin{equation}\label{nondegen}
[X_{\ell,j},X_{m,j}]\equiv a_{\ell
m}^jT
\qquad 
\mod (X_{\ell,j})
\end{equation}
then the matrix $(a_{\ell m}^j)$ is non-degenerate. 
\end{itemize}
\begin{remark} It is a result that goes back to Tanaka
(\cite{Tanaka}) and used to advantage in the work of the present
authors in many places that under the non-degeneracy assumption
(\ref{nondegen}), one {\it always} may modify the
given vector field $T$ by adding multiples of the $X_{j,\l}$ in
such a way that 
\begin{equation}\label{Tcommutes}
[X_{j,\l}, T] \equiv 0 \mod (X_{1,\l}, \ldots , X_{1,\l}).
\end{equation}
\end{remark}

\begin{definition} We call such a family
$(X_{\ell,j}, T)$ of systems of vector fields a
compatible family. 
\end{definition}

Now, we may assume that each $(V_j)$ is the domain
of a coordinate chart. So in each $V_j$ and for
every $s\geq 0,$ we may consider an elliptic
pseudodifferential operator of order $s$ which we
denote by $\Lambda_j^s.$

Let us fix a family $(\phi_j)_{j=1,\ldots p}$ such
that 
\begin{equation} \phi_j \in \D(V_j), \qquad 0\leq
\phi_j \leq 1, \qquad \sum \phi_j\equiv 1 \hbox{ on
} M.
\end{equation}

Let $(\rho_j)_{j=1,\ldots p}$ be another family
such that such that 
\begin{equation} \rho_j \in \D(V_j), \qquad 0\leq
\rho_j \leq 1, \qquad \rho_j\equiv 1 \hbox{ on
supp} \;\phi_j.
\end{equation}

Now one has, say for $t\geq s\geq 0,$
\begin{equation}\label{t,t-s}
\|\phi_jv\|_t \lesssim (\|\rho_j\Lambda_j^s \phi_j
v\|_{t-s} + \|\phi_jv\|_0), \forall v \in C^\infty (M)
\end{equation}
where $\|\;\;\|_t$ denotes the Sobolev norm. 

So, now, one has 
\begin{equation}\label{t,t-s,bis}
\|v\|_t\lesssim \sum_j\|\phi_jv\|_t \lesssim
\sum_j (\|\rho_j\Lambda_j^s
\phi_j v\|_{t-s} + \|\phi_jv\|), \forall v \in
C^\infty (M)
\end{equation}

Now, in each $V_j,$ we consider the operator
considered in the introduction (depending on the
given $u\in C^\infty (M)$): 
\begin{equation}
\begin{cases}
P_j  = \sum_{\ell = 1}^{2n} (Y_{\ell,j}^2 +
a_{\ell,j}Y_{\ell,j}+b_j)\\\\
a_{\ell,j}, b_j \in \A(V_j)
\end{cases}
\end{equation}
and assume that 
\begin{equation}
P_j u \in \A(V_j), \quad \forall j.
\end{equation}

Finally let us denote by $(\quad , \quad )_s$ the
$s-$Sobolev scalar product (in each $V_j,$ when one has
functions with compact support in $V_j$) and remember
the following:
\begin{equation}\label{pm1/2}
\forall \delta >0, \exists C_\delta:\forall w\in
C_0^\infty,  
\|w\|_s^2 \leq \delta \|w\|_{s+1/2}^2 + C_\delta \|w\|^2_0.
\end{equation}

\section{Maximal Estimates}
\renewcommand{\theequation}{\thesection.\arabic{equation}}
\setcounter{equation}{2}
\setcounter{theorem}{0}
\setcounter{proposition}{0}  
\setcounter{lemma}{0}
\setcounter{corollary}{0} 
\setcounter{definition}{0}
\setcounter{remark}{1}

Our aim in this section is to prove the following: 
\begin{theorem}
We have the following maximal estimates for $s\geq 0:$
\begin{equation}\tag{$3.1_s$}\label{3.1s}
    \|v\|_{{s+1/2}}^{2} + 
    \sum_{j,\ell}\|X_{j,\ell}\phi_{\ell}v\|_{s}^{2} 
\lesssim 
    \sum_{\ell}|(\phi_{\ell} P_{\ell}v, \phi_{\ell}v)| +
\|v\|^{2}_0, 
\end{equation}
for $v\in C^{\infty} (M)$ and
\begin{equation}\tag{$3.2_s$}\label{3.2s}
    \|v\|_{{s+1}}^{2} + 
    \sum_{j,\ell}\|X_{j,\ell}\phi_{\ell}v\|_{s+1/2}^{2}
   \lesssim 
    \sum_{\ell}\|\phi_{\ell} P_{\ell}v\|_{s}^{2} + \|v\|^{2}_0 
\end{equation}
%\end{theorem}
%\end{document} 
for $v\in C^{\infty}(M)$ and in fact
\begin{multline}\label{3.3s}\tag{$3.3_s$}
 \|v\|_{s+1}^{2} + 
    \sum_{j,\ell}\|X_{j,\ell}\phi_{\ell}v\|_{s+1/2}^{2}
+\sum_{j,k,\ell}\|X_{j,\ell}X_{k,\ell}\phi_{\ell}v\|_{s}^{2}\\
\lesssim 
    \sum_{\ell}\|\phi_{\ell} P_{\ell}v\|_{s}^{2} + \|v\|^{2}_0 
\end{multline}
for $v\in C^{\infty}(M).$
\end{theorem}
\begin{proof}
    The proof is known when written for functions with compact 
    support in coordinate charts. This is a global version. Let us 
    first show the statements at level $s=0.$ For simplicity we 
    take $\ell$ fixed and set $X_{j,\ell}= X_{j}, j=1, \ldots 2n$ and 
    $\phi = \phi_{\ell}.$ We have
\[
\sum_j \|X_j\phi v\|^2_0 = \sum (X_j\phi v,
X_j\phi v) \lesssim \sum (Y_j\phi v,
Y_j\phi v)
\]
because $H$ is invertible. Note that $\lessim$ may indicate a
constant which depends on $u$ and its first few derivatives. Now 
\begin{multline*}
(Y_j\phi v,Y_j\phi v) = -(Y_j^2 \phi v, \phi v) +
(\theta_j \phi v, \phi v), \qquad \theta_j \in C^\infty
(V_j)\\
=-([Y_j^2,\phi]v,\phi v) + (\phi Y_j^2 v, \phi v) +
\mathcal{O} (\|\phi v\|_0\|Y_j\phi v\|_0). 
\end{multline*}

Now, using $[Y_j^2,\phi ] = 2Y_j[Y_j, \phi ] - [Y_j,
[Y_j,\phi ]]$
we easily obtain
\begin{equation}\label{[Ysquared,phi]}
|([Y_j^2,\phi ] v, \phi v)| \lesssim C _1 \|v\|^2_0 +
{1\over C_0}\sum_j \|X_j \phi v\|^2_0.
\end{equation}
Thus, 
\[
\sum_j \|X_j \phi v\|^2_0 \lesssim |(\phi Pv, \phi v)| +
C_1\|v\|^2_0 + {1\over C_0}\sum_j \|X_j \phi v\|^2_0.
\]

Now use 
\begin{equation}
\|\phi v\|_{1/2} \lesssim \sum_j\|X_j\phi v\|^2_0 + \|\phi
v\|^2_0 \qquad \forall v \in C^\infty (M) 
\end{equation}
(see J. J. Kohn \cite{K}) to obtain (\ref{3.1s}) in
case $s=0.$ 

Now we can deduce $(3.2_0)$ from $(3.1_0)$ in
the following way:  using (\ref{t,t-s}) and
(\ref{t,t-s,bis}), we have
\begin{multline}
\sum_\ell\|\phi_\l v\|_1^2 + \sum_{j,\l}\|X_{j,\l}\phi_\l
v\|_{1/2}^2 \lessim \sum_\l
\|\rho_\l\Lambda_\l^{1/2}\phi_\l v\|_{1/2}^2\\
+\sum_{j,\l}\|\rho_\l\Lambda_\l^{1/2}X_{j,\l}\phi_\l v\|^2_0
+
 \|v\|^2_0 + \sum_{j,\l}\|X_{j,\l}\phi_\l v\|_0^2\\
\lessim \sum_\l \|\rho_\l\Lambda_\l^{1/2}\phi_\l
v\|^2_{1/2} + \sum_{j,\l} \|X_{j,\l}\rho_\l\Lambda_\l^{1/2}\phi_\l
v\|^2_0 + \sum _\l \|\phi_\l v\|^2_{1/2}+\sum_\l\|X_{j,\l}
\phi_\l v\|^2_0.
\end{multline}
The last two sums are easy to handle. The first two sums
are (from the first part of the theorem at level $s=0$),
less than 
\begin{equation}
\sum_\l |(\rho_\l P_\l \Lambda _\l^{1/2} \phi_\l v,
\rho_\l
\Lambda _\l^{1/2}\phi_l v)| + \sum_\l \|\phi_l v\|^2_{1/2}
\end{equation}
Now, for simplicity we forget the index $\l$ and consider 
\begin{equation}\label{a}
\rho P \Lambda ^{1/2} \phi v = [\rho P, \Lambda^{1/2}
\phi] v + \Lambda^{1/2}\phi P v
\end{equation}

Now as we obtained (\ref{[Ysquared,phi]}) we have 
\begin{equation}
|([\rho P, \Lambda^{1/2}\phi]v, \rho \Lambda ^{1/2} \phi
v)| \leq C_1\|v\|^2_{1/2} + {1 \over C_0}\sum\|X_j \phi
v\|^2_{1/2}\label{[rhoP,Lambdaphi]}
\end{equation}
By taking $C_0$ large enough and using (\ref{pm1/2}) we
have the desired inequality, because the term
$|(\Lambda^{1/2}\phi Pv, \rho\Lambda^{1/2}\phi v)|$ is
less than 
$C_1\|\phi Pv\|^2 + {1 \over C_0}\sum\|\phi v\|^2_1$
(with $C_0$ large, $C_1$ depending on $C_0$ as usual). 

This proves (\ref{3.2s}) with $s=0.$ To bound also the third 
term on the left in (\ref{3.3s}) for $s=0,$ we argue as follows:
first the function $X_{j,\l}v$ is inserted in place of $v$ in
(\ref{3.1s}) with $s=0$  and an error of the type
$C\sum_\l \|v\|_1^2$ is introduced through  a bracket of the form
$([X,X]v, X^2v).$ While this is a new error, it is already
controlled by (\ref{3.2s})s, which completes the proof for $s=0.$ 

Our aim is to prove (\ref{3.1s}) and deduce
(\ref{3.2s}) from (\ref{3.1s}) as before. 

First observe that (in view of (\ref{t,t-s}))
\[
\|v\|^2_{s+{1/2}}\lessim
\sum_\l\|\rho_\l\Lambda_\l^s\phi_\l v\|^2_{1/2} +
\|v\|^2_0 \qquad \hbox{ and }
\]
\[
\sum_\l \|X_{j,\l}\phi_\l v\|^2_s\lessim
\sum_\l\|X_{j,\l}\rho_\l\Lambda_\l^s\phi_\l v
\|^2_0 + \|v\|^2_0 \lessim
\]
\[
\lessim
\sum_\l\|X_{j,\l}\rho_\l\Lambda_\l^s\phi_\l v
\|^2_0 +
\sum_\l\|\rho_\l\Lambda_\l^s\phi_\l v\|_{1/2}^2 + \|v\|^2_0.
\]

Now, we use ($3.1_0$) to obtain
\begin{equation}\label{vs+1/2}
\|v\|^2_{s+1/2} + \sum_\l\|X_{j,\l}\phi_\l v\|_s^2\lessim
\sum_\l |(\rho_\l P_\l\Lambda_\l^s\phi_\l v,
\rho_\l\Lambda_\l^s \phi_\l v)|.
\end{equation}

In view of (\ref{pm1/2}), the term $\|v\|^2_s$ may be
replaced by $\|v\|^2.$

Now we consider the first term in the second member of
(\ref{vs+1/2}) and write: $\rho_\l P_\l \Lambda_\l^s
\phi_\l v = \rho_\l [P_\l , \Lambda_\l^s \phi_\l ]v +
\rho_\l \Lambda_\l^s \phi_l P_\l v.$

So one sees that one is reduced to study $(\rho_\l [P_\l, 
\Lambda_\l^s \phi_\l ]v, \rho_\l \Lambda_\l^s \phi_l
 v)$, because one has easily 
\begin{multline}
|(\rho_\l 
\Lambda_\l^s \phi_\l P_\l v, \rho_\l \Lambda_\l^s \phi_l
 v)| \leq C|(\phi_\l P_\l v, \phi_\l v)_s|\\
+\delta\left\{\sum_\l \|\phi_l v\|^2_{s+1/2} + \sum_\l
\|X_{j,\l }\phi_l v\|^2_{s}\right\} + C_\delta \|v\|^2_0.
\end{multline}

Now again forget the index $\l$ and consider 
\[
[P,\Lambda^s\phi]=\sum [Y_j^2,\Lambda^s\phi ] = \sum 2Y_j
[Y_j,\Lambda^s\phi] - [Y_j,[Y_j,\Lambda^s\phi]]
\]

Then one has, as in
(\ref{[Ysquared,phi]}), 
\[
|(\rho_\l [P_\l ,\Lambda_\l^2\phi_\l]v, \rho_\l
\Lambda_\l^s\phi_\l v)| \leq {1\over C_0}\sum_j\|X_{j,\l
}\rho_\l \Lambda_\l^s \phi_\l v\|^2_0 +C_1\|\rho_\l
\Lambda_\l^s\phi_\l v\|^2_0
\]
where $C_1$ depends on $C_0,$ as usual. 

Then, again using (\ref{pm1/2}) and taking $C_0$ big
enough, the first member of (\ref{vs+1/2}) is less than
$C(\sum_\l |(\phi_\l P_\l v, \phi_\l v)_s| + \|v\|^2_0)$
for some $C>0.$

Now, we want to prove (\ref{3.2s}) using (\ref{3.1s}) as
we did for $s=0.$

One has, as in that case, 
\[\sum_\l \|\phi_\l v\|^2_{s+1} +
\sum_{j,\l}\|X_{j,\l}\phi_\l v\|^2_{s+1/2} \lessim
\sum_\l \|\rho_\l \Lambda_\l^{1/2}\phi_\l v\|^2_{s+1/2}
\]
\[+\sum_{j,\l}\|\rho_\l \Lambda_\l^{1/2}X_{j,\l}\phi_\l
v\|^2_{s} + \|v\|^2_{s+1/2} +
\sum_{j,\l}\|X_{j,\l}\phi_\l v\|^2_s
\]
\[\lessim \sum_\l \rho_\l \Lambda_\l^{1/2}\phi_\l
v\|^2_{s+1/2} + \sum_{j,\l}\|X_{j,\l}\rho_\l
\Lambda_\l^{1/2}\phi_\l v\|^2_s + \sum_\l \|\phi_\l
v\|^2_s
\]
\[ + \sum_{j,\l}\|X_{j,\l}\phi_\l v\|^2_s.
\]
The worst terms are the first two ones, from
(\ref{3.1s}). They are less than
\[\sum_\l |(\rho_\l P_\l \Lambda_\l^{1/2}\phi_\l v,
\rho_\l \Lambda_\l^{1/2}\phi_\l v)_s| + \|v\|^2_0
\] 
Now the end of the proof of (\ref{3.2s}) follows the
lines of the end of the proof in the case $s=0,$ and 
the proof of (\ref{3.3s}) is also as before.
\end{proof}    
\begin{remark} This proof of the global version only
requires careful computations.
\end{remark}
\begin{corollary} Let $T$ be a global, real analytic,
non-zero vector field on $M$ complementary to the $X$ and
satisfying (\ref{nondegen}). Then 
\[\|Tv\|_s^2 \lessim \sum_\l \|\phi_\l P_\l v\|^2_s +
\|v\|^2_0 \qquad \forall v \in C^\infty(M)\]
\end{corollary}
\begin{remark}
The existence of such a global $T$ has been shown
(\cite{Tartakoff1978}, \cite{Tartakoff1981}) when $M$ is a compact
real compact CR manifold.
\end{remark}

\section{High Powers of the Vector Field $T.$}
\renewcommand{\theequation}{\thesection.\arabic{equation}}
\setcounter{equation}{0}
\setcounter{theorem}{0}
\setcounter{proposition}{0}  
\setcounter{lemma}{0}
\setcounter{corollary}{0} 
\setcounter{definition}{0}
\setcounter{remark}{1}

The overall strategy is to use the maximal estimates above with
$v$ replaced by $T^p u.$ Once one has control over high $T$
derivatives of the solution the other derivatives follow by
standard techniques.  

Now the vector field $T$ being global, dealing with
just a bounded number of vector fields $X_j$ that are only locally
defined is not a very delicate issue. For instance, the above
Corollary may be strengthened to 

\begin{equation}
\|Tv\|_s^2+\|v\|_{s+1}^{2} + 
\sum_{j,\ell}\|X_{j,\ell}\phi_{\ell}v\|_{s+1/2}^{2}
+\sum_{j,k,\ell}\|X_{j,\ell}X_{k,\ell}\phi_{\ell}v\|_{s}^{2}
\end{equation}
\begin{equation*}
\lessim \sum_\l \|\phi_\l P_\l v\|^2_s +
\|v\|^2_0 \qquad \forall v \in C^\infty(M)
\end{equation*}
or, in the form we will use it, for any integer $r,$
\begin{equation}\label{apriori:TrHs}
\|T^{r+1}v\|_s^2 + 
\sum_{j,\ell}\|X_{j,\ell}\phi_{\ell}T^rv\|_{s+1/2}^{2}
+\sum_{j,k,\ell}\|X_{j,\ell}X_{k,\ell}\phi_{\ell}T^rv\|_{s}^{2}
\end{equation}
\begin{equation*}
\lessim \sum_\l \|\phi_\l P_\l T^rv\|^2_s +
\|T^rv\|^2_0 \qquad \forall v \in C^\infty(M)
\end{equation*}

\begin{proposition} There exist constants $C, C_u, C_h$ such
that for all $r,$ 
\begin{equation}\label{pureTpowers}
{1\over r!}\{\|T^{r+1}v\|_s^2 +
\|X^2T^rv\|_{s}^{2}\} \leq 4^rC^r C_u^r C_h^r
\end{equation}
\end{proposition}

\begin{proof}

In view of (\ref{apriori:TrHs}), our (only) task is to
commute
$T^r$ past $P_\l$ with errors which can be recursively estimated
to grow `analytically' since then after specializing $v$ to $u,$ we
know that $P_u u \in
\A(M).$

Now we have the crucial relationship (in $V_j$)
$$[T,X_{j,\l}] = \sum_k c_{jk}X_{k,\l}$$
with $c_{jk}$ real analytic functions. Since $P_\l$ is a quadratic
polynomial in the $X_{j,\l},$ with coefficients $h(x,u)$ which are
real analytic functions of the spatial variables $x$ {\it and} the
solution $u(x),$ we may write 
$$[T,P_\l] = \sum [T, h(x,u)X_{k,\l}X_{j,\l}] = (Th)X^2 +
h[T,X^2]$$
$$=(Th)X^2 + h \{\tilde{a}X^2 + {\tilde{\tilde{a}}}X +
{\tilde{\tilde{\tilde{a}}}})$$
with analytic functions $\tilde{a}, \tilde{\tilde{a}}$ and 
$\tilde{\tilde{\tilde{a}}}.$ Generically, then, we may write 
\begin{equation}\label{genericTPell}
[T,P_\l] = \sum ((T)h) X^2
\end{equation}
where $h((T)h)$ denotes (at most) a first derivative (in $(x,t)$)
 of $h(x,u(x))$ times one of a finite collection of
analytic functions of $x,$ namely the coefficients of the $X$
in the bracket $[T,X]$ mentioned above in (\ref{nondegen}) (and
possibly one derivative of this coefficient). There may also be
fewer than two
$X's$ on the right in (\ref{genericTPell}). For the rest of the
paper we will assume for simplicity that $h=h(u).$

Next, we will need an expression for the more complicated bracket 
$$[T^r,
P_\l]=\sum_{r'=0}^{r-1}T^{r'}[T,P_\l]T^{r-r'-1}
=\sum_{r'=0}^{r-1}T^{r'}\circ h'X^2\,T^{r-r'-1}
$$
only we move the $h'$ to the very left yet leave $X^2$ for the
moment wherever they are, which we denote by enclosing the $X^2$
in parentheses and placing them on the left. That is, 
\begin{equation}
\label{bracketTrPell}
[T^r,P_\l]=\sum_{r'=1}^{r}{r\choose
r'}(T^{r'}h)(u(x))(X^2)T^{r-r'}.
\end{equation} 
or
$$
{[T^r,P_\l]\over r!}=\sum_{r'=1}^{r}{(T^{r'}h)(u(x))\over r'!}
\;{(X^2)T^{r-r'}\over (r-r')!}.$$

Now for the term $(T^{r'}h)(u(x)),$ we will
use the Fa\`a di Bruno formula or rather crude bounds for the
results: writing
$$D^k_x g(u(x)) = (u'D_u+D_x)^{k-1}u'D_u g(u(x));$$
writing this crudely as 
$$D^k_x g(u(x)) = ((u'\sigma +D_x)^{k-1}u'_{|\sigma = D_g})g,$$
i.e., $\sigma$ becomes a `counter' for the number of derivatives
received by $g.$ Then this is {\it at worst} 
\begin{equation}\label{approxFdB}\sum_{k'}{k \choose
k'}g^{(k-k')}(D_x^{k'}{u'}^{k-k'}).\end{equation}  
Finally, we must
analyze expressions such as $D^a {u'}^b:$
$$ D^a {u'}^b = \sum_{a_1+\ldots + a_b=a}
{a \choose a_1, \ldots a_b} {u'}^{(a_1)}\ldots {u'}^{(a_b)}$$
where ${a \choose a_1, \cdots a_b}$ denotes the multinomial
expression 
$${a \choose a_1, \ldots a_b} = {a! \over a_1!\ldots a_b!(a-\sum
a_j)!}= {a! \over a_1!\ldots a_b!}$$
since the $\sum a_j = a.$ We have 
$$ {D^a {u'}^b \over a!} = \sum_{a_1+\cdots + a_b=a} 
{{u'}^{(a_1)}\over a_1!}\cdots {{u'}^{(a_b)}\over a_b!}$$
Thus we have 
\begin{equation}\label{bracketwithD^r}
{[T^r,P]\over r!}=\sum_{r'=1}^{r}
{(T_{x,t}^{r'}h(u(x,t)))\over r'!}
\;{(X^2)T^{r-r'}\over (r-r')!}.
\end{equation}
and so (cf. (\ref{approxFdB})):
$${(T^{r'}) h(u(x))\over r'!}
\sim
\sum_{r'-r''\geq 1}^r{h^{(r'-r'')}\over
{(r'-r'')}!}\;{(T^{r''}({u'}^{r'-r''}))\over r''!}$$ 
$$=\sum_{r'-r''\geq 1}^r{h^{(r'-r'')}\over
{(r'-r'')}!}\sum_{\sum_1^{r'-r''} r''_j = r''} 
{T^{r''_1}{u'}\over r''_1!}\cdots {T^{r''_{r'-r''}}{u'}\over
r''_{r'-r''}!}$$ or in all, 
\begin{equation}
\label{[Tr,P]}
{[T^r, P]v\over r!}
=\sum_{{{r\geq r'-r''\geq 1}\atop \sum_{j=1}^{r'-r''} r''_j =
r''}\atop {(\sum_{j=1}^{r'-r''} (r''_j+1) =
r'})}{h^{(r'-r'')}\over
{(r'-r'')}!}{T^{r''_1}{u'}\over r''_1!}\cdots
{T^{r''_{r'-r''}}{u'}\over r''_{r'-r''}!}\;{(X^2)T^{r-r'}v\over
(r-r')!}.
\end{equation}
To simplify the argument we have dropped the localizing
functions, since for global arguments when these functions always
appear on the left of the norms they may be brought out easily and
replaced by another partition of unity, with new X's if needed;
also have ignored the order of the $X's$ and $T's,$ indicating
this by putting the $X^2$ in parentheses, not to indicate that
they may not be present (though they may not) but that they may
appear with several $T's$ to the left and more to the right. We
have also dropped all subscripts. Schematically, we may then write
(\ref{apriori:TrHs}) together with (\ref{[Tr,P]}) as

\begin{multline}
\label{}
\|T^{r+1}v\|_s + \|X^2T^rv\|_{s}
\lessim \|P T^rv\|_s +
\|T^rv\|_s \leq\\
 \lesssim \|T^r P v\|_s + \|T^rv\|^2_0 + \|[T^r,P]v\|_s 
\qquad \forall v \in C^\infty(M) 
\end{multline}
and so
\begin{multline}
\label{normofproduct}
{1\over r!}\{\|T^{r+1}v\|_s +
\|X^2T^rv\|_{s}\}
\lesssim {1\over r!}\{\|T^r P v\|_s +
\|T^rv\|_s\} \;+ \\
+\sum_{{{r\geq r'-r''\geq 1}\atop \sum_{j=1}^{r'-r''} r''_j =
r''}\atop {(\sum_{j=1}^{r'-r''} (r''_j+1) =
r'})}\left\|{h^{(r'-r'')}\over {(r'-r'')}!}\;{T^{r''_1}{u'}\over
r''_1!}\cdots {T^{r''_{r'-r''}}{u'}\over
r''_{r'-r''}!}\;{(X^2)T^{r-r'}v\over (r-r')!}\right\|_{s}.
\end{multline}

Now for our value of $s, H^s$ is an algebra, and so
the norm of the product of derivatives of copies of $u$ may be
replaced by the product of the norms, each of which will have the
form of one of the terms on the left hand side of
(\ref{normofproduct}). 

Specializing to $v=u$ and bounding the norm of the product by the
product of the norms we observe that except for the term
involving derivatives of $h,$ all other terms are of the same form
since one $T$ derivative and two $X$ derivatives carry the same
weight on the left hand side of (\ref{normofproduct}):

\begin{multline}
\label{productofnorms}
{1\over r!}\{\|T^{r+1}v\|_s +
\|X^2T^rv\|_{s}\}
\lesssim {1\over r!}\{\|T^r P v\|_s +
\|T^rv\|_s\} \;+ \\
+\sum_{r\geq r'-r''\geq 1}C^{r'-r''+2}\left\|{h^{(r'-r'')}\over
{(r'-r'')}!}\right\|_s \left\|{(X^2)T^{r-r'}v\over
(r-r')!}\right\|_{s}\prod_{\sum_{j=1}^{r'-r''}
r''_j = r''}\left\|{T^{r''_j}{u'}\over
r''_j!}\right\|_s.
\end{multline}
The constant includes a power of $C$ for each norm that
follows. Note that since the derivatives on $h$ are of that order,
this constant will be included with the analyticity constant for
$h,$ and in the future constants with exponents comparable to
the number of derivatives on a function known to be analytic
will be permitted without comment. 

Note that the term in the product with $(X^2)$ 
is analogous to the extra $T$ derivative on each of the other
terms. Hence these terms are similar to the left hand side, and
could be handled at once inductively except for counting the
number of them, but it is simpler to
iterate (\ref{productofnorms}) directly, at least until all terms
have order less than $r/2.$

Since there can be at most one term of order larger than $r/2,$ 
after the next `pass' we observe that the product will look just
like the right hand side of (\ref{normofproduct}) again, except
that there will be one more norm of derivatives of $h.$ 
 
That is, applying (\ref{productofnorms}),
with $r$ replaced by $r-r',$ to the  
term in (\ref{productofnorms}) with $(X^2),$ we have 
\begin{multline}
\label{productofnorms2}
{1\over (r-r')!}\{\|T^{r-r'+1}v\|_s +
\|X^2T^{r-r'}v\|_{s}\}\lesssim \\
\lesssim {1\over {r-r'}!}\{\|T^{r-r'} P v\|_s +
\|T^{r-r'}v\|_0\} \;+ \\
+\sum_{r-r'\geq \rho'-\rho''\geq 1}C^{\rho'-\rho''+2}\left\|{h^{(\rho'-\rho'')}\over
{(\rho'-\rho'')}!}\right\|_s\left\|{(X^2)T^{r-r'-\rho'}v\over
(r-r'-\rho')!}\right\|_{s}
\prod_{\sum_{j=1}^{\rho'-\rho''} \rho''_j = \rho''}
\left\|{T^{\rho''_j}{u'}\over
\rho''_j!}\right\|_s
\end{multline}
we find $r\geq r'+\rho'-r''-\rho''$ and
$\sum_{j=1}^{r'-r''}\sum_{k=1}^{\rho'-\rho''}
(r_j^{''}+\rho_k^{''})= r''+\rho''$
so if we set $s'=r'+\rho'$ and $s''=r''+\rho'',$ we have
a sum over $s'-s''$ and $\sum_{j+k=2}^{s'-s''} s_{j+k}''=s''$
subject to the obvious subdivisions.

That is, over $r'+\rho'=s', r''+\rho''=s'',$
\begin{multline}
{1\over r!}\{\|T^{r+1}v\|_s +
\|X^2T^rv\|_{s}\}
\lesssim \\
\lesssim \sum_{r'\geq 0}{1\over (r-r')!}\{\|T^{r-r'} P v\|_s +
\|T^{r-r'}v\|_s\} \;+ \\
\label{productofnorms3}+\sum_{s'=r'+\rho', s''=r''+\rho''
\atop {r\geq r'-r''\geq 1\atop {r-r'\geq \rho'-\rho''\geq 1\atop
{\sum_{j=1}^{s'-s''} s''_j =
s'' \atop {(\sum_{j=1}^{s'-s''} (s''_j+1) =
s'})}}}}\left\|{C^{r'-r''+2}h^{(r'-r'')}\over
{(r'-r'')}!}\right\|_s\left\|{C^{\rho'-\rho''+2}h^{(\rho'-\rho'')}\over
{(\rho'-\rho'')}!}\right\|_s\,\times \\
\times \left\|{T^{s''_1}{u'}\over
s''_1!}\right\|_s\cdots
\left\|{T^{s''_{s'-s''}}{u'}\over
s''_{s'-s''}!}\right\|_s\left\|{(X^2)T^{r-s'}v\over
(r-s')!}\right\|_{s}.
\end{multline}

Note that in using the fact that $H^s$ is an algebra, i.e.,
$\|fg\|_s\leq B\|f\|_s\|g\|_s,$ we have absorbed the algebra
constant with the constant $C$ inside the norms of $h.$ We
further estimate the norms of derivatives of $h$ (noting that
each occurrence contains at least one such derivative) by
\begin{equation}\label{eachsummand}\|C^{\ell+2}h^{(\ell)}(x,y,u)\|_s\leq
C_h^\ell \,\ell!\end{equation}

We are nearly ready to iterate this procedure until even the last
term has order less than $r/2;$ for except for the sum (the
number of terms), each term has a bound which is stable in the
number of iterations, namely the last right hand side above is
bounded by 
\begin{equation}
\label{product}
\sum C_h^{t}\prod
\left\{\left\|{T^ku'\over k!}\right\|_s
\hbox{ or } \left\|{X^2T^ku\over k!}\right\|_s\right\}
\end{equation} 
where the
sum of the $k+1$ is at most $r$ and $t\leq r$ is the number of
terms in the product. 

As for the sum, whether after a single full pass or multiple
ones, the number of terms corresponds {\it at most} to the number
of ways to partition $r$ derivatives among at most
$r$ functions, generally many fewer. Denoting by $D$ a
derivative ($r$ of them) and by $u$ a copy of $u$ ($t$ of
them) we are faced with the number of ways to `identify' or select
$t$ items (the $u's$) from among $r+t$ items (the $D$'s and
$u$'s) with the understanding that in an expression such as 
\begin{equation}\label{partitions}
\underbrace{\underbrace{DDDDD}_{r_1}u
\underbrace{DDDDD}_{r_2}u
\underbrace{DDDDD}_{r_3}u
\cdots
\underbrace{DDDDD}_{r_t} u}_{r \; D's \;{\rm and }\; t(\leq
r)\, 
\;u's}
\end{equation} 
the $D$'s differentiate only the first $u$ which follows. The
answer is that there are certainly not more than ${r+t \choose
t}\leq 2^{r+t}\leq 2^{2r}=4^r$ ways.  

Finally, since we may iterate this procedure until the maximal
order of differentiation on $u$ is $1$ or $2,$ and bound this
small number of derivatives by a constant (with at most $r$ such
terms, naturally - that's all the derivatives there were). Thus the
left hand side of (\ref{normofproduct}) is bounded by: 
\begin{equation}\label{pureT}
{1\over r!}\{\|T^{r+1}v\|_s^2 +
\|X^2T^rv\|_{s}^{2}\} \leq 4^rC^r C_u^r C_h^r
\end{equation}
which clearly yields analytic growth (of $T$ derivatives) of the
solution $u$ since $C_u$ depends only on the first few
derivatives of
$u$ and $s$ is taken just large enough to ensure that $H^s$ is an
algebra. 
\end{proof}

\section{Mixed Derivatives - the case of global X}
\renewcommand{\theequation}{\thesection.\arabic{equation}}
\setcounter{equation}{0}
\setcounter{theorem}{0}
\setcounter{proposition}{0}  
\setcounter{lemma}{0}
\setcounter{corollary}{0} 
\setcounter{definition}{0}
\setcounter{remark}{1}

To finish the proof in the case where the vector field(s) $X$ are
globally defined it remains to show that we may estimate mixed
derivatives as effectively as we did the high
$T$ derivatives. The result of Helffer and Mattera show
that it {\it would} suffice to handle pure powers of the
vector field
$X,$ but mixed derivatives will invariably enter through
brackets of pairs of $X$'s. Thus this we start by using the {\it a
priori} estimate (\ref{3.1s})-(\ref{3.3s}) with
$v$ replaced by
$\phi X^r$ (and later by a mixture of derivatives in $X$ and in
$T$). What ultimately happens is that brackets of pairs $X's$ will
produce
$T's$, but at most half as many, and we will be led back to
(nearly pure) powers of
$T.$ The
non-linearity of the problem introduces nothing new in this
overall pattern. 

When the $X$'s are {\em globally defined}, for example in $\C
^2,$the powers of
$X$ are treated like powers of $T$ (e.g., with respect to the use
of the F\`aa di Bruno formula, especially) with the one change
that an additional type of term will appear: starting with
$X^{r+2}v$ there will appear as an error $\underline{r}$ copies
of $X^rTv$ when two $X$'s bracket to give a $T.$ And this
effect, the only new feature, will be repeated until all or nearly all
the $X$'s are exhausted. That is, we have the new scheme 
$$X^r\rightarrow C^{r/2} r!!(X^2)T^{r/2},$$
where we recall the definition $r!!=r(r-2)(r-4)\ldots \sim
C^{r/2}(r/2)!$
But this is not a problem, since we have just treated essentially
pure powers of $T$ above in (\ref{pureT}). 

Rather than write this case out in more detail, we proceed to the
next section where the {\it problem} is global but the vector
fields
$X$ are only locally defined. This case incorporates many of
the features of a fully local proof, though fortunately not
all! Note that the
$T$ vector field is still required to be globally given. 

\section{Mixed Derivatives - when $X_j$ are only locally
defined}
\renewcommand{\theequation}{\thesection.\arabic{equation}}
\setcounter{equation}{0}
\setcounter{theorem}{0}
\setcounter{proposition}{0}  
\setcounter{lemma}{0}
\setcounter{corollary}{0} 
\setcounter{definition}{0}
\setcounter{remark}{1}

When the vector fields $X$ are only locally defined, we cannot
afford to change freely from one coordinate patch to another and to
another basis of
$X's$ each time a localizing function arising from a partition of
unity is differentiated, since the constants counting the number
of terms and the coefficients would grow far too fast, namely
roughly $C^r$ at {\em each} step.  we will
need a suitable localization of high powers of the
$X.$ While one might suspect otherwise, we will 

We will thus work in a single coordinate patch, drop all
subscripts $\l$ (and $j$ and $k,$ for simplicity), and in place
of $v$ in the {\it a priori} estimates substitute
$\Psi X^ru,$ where the localizing function $\Psi$
will be specified further below. It will not need to be
differentiated very often, but the band in which it goes from
being identically equal to one to being identically zero will be
of a precise width, as will subsequent localizing functions which
will be introduced below. The general result on families of
localizing functions is given by a result of Ehrenpreis
(\cite{Ehrenpreis1960}):

\begin{proposition}
For any two open sets $\Omega_0\Subset\Omega_1,$ with
separation $d={\hbox{dist.}}(\Omega_0,\Omega_1^{c})$ and any
natural number $N,$ there exists a universal constant
$C$ depending only on the dimension and a function
$\Psi=\Psi_{\Omega_0,\Omega_1,N}\in C_0^\infty (\Omega_1),
\Psi\equiv 1 \hbox{ on }\Omega_0$ 
with 
\begin{equation}\label{fns:Ehrenpreis}
|D^\beta \Psi| \leq \left({C\over
d}\right)^{|\beta|+1}N^{|\beta|},
\qquad |\beta| \leq 2N,
\end{equation}
\end{proposition}

\noindent though in this paper we will take
$N=4$ at most; thus the analyticity to be shown in $U_0$ will be
reduced to combining the bounds on $\|T^{r+1}u\|_s$ obtained
in a previous section with the bounds 
$$\sum\| X^2 \Psi X^r u\|_s + \|T\Psi
X^r u\|_s \leq C_\Omega^{r+1}r! $$
with $\Psi\equiv 1$ on the set where we want to prove
analyticity.

To do this, we start with the {\it a priori} estimates as before: for
$v$ of compact support where the $X_j$ are defined, and any fixed
$s,$ we have (\ref{3.3s}) in the form: 
\begin{equation}\label{est:apriori}
 \|X^2 v\|_{s}^{2} + 
    \|X v\|_{s-1/2}^{2}
+\|v\|_{s-1}^{2}
\lesssim 
   \| Pv\|_{s}^{2} + \|v\|^{2}_0. 
\end{equation}

Note that we have dropped all subscripts but are working with several
$X$'s. Naturally we could add a term $\|Tv\|_{s}$ to the left hand
side using the non-vanishing of the Levi form but it will not help us
here as it did above in handling high powers of $T$
applied to the solution $u.$

This estimate will be applied to $v=\Psi X^r u$ and then on the right
we will write $P\Psi X^r u$ in terms of $\Psi X^r Pu$ modulo an
error, 
namely the commutator of $aX^2$ with $\Psi X^r,$ suitably expanded. 

Now the crucial brackets, analogous to 
(\ref{bracketTrPell}), will be written
\begin{equation}\label{exp:[P, Psi X^r]-1}
[ P, \Psi X^r]v = a_u\, [X^2,\Psi ] X^rv + a_u\,\Psi
[X^2,X^r]v + \Psi [a_u\,, X^r] X^2v,
\end{equation}
where coefficients depending on the
solution
$u$ (those arising in
$P$ and here denoted $a_u$) are subscripted with
$u$ while those which depend only on the spatial
variables are not subscripted. Now 
$$ [X^2,\Psi] = 2\Psi' X+ 
 \Psi'',$$
(and notice that at most two derivatives appear on $\Psi$ and that
these will fall to the left of all other terms in the bracket and will
be changed with each iteration) and 
\begin{equation}\label{[X,X]=aT}
[X,X]=aT \qquad {\hbox{ and so }}
\end{equation}
$$
[X^2,X^r] =
\underline{Cr}{\hbox{ terms }}\;[X,X]X^r+\cdots $$
independent of $u$.
Here underlining a coefficient indicates the number of
terms of the given type which occur and the
$\cdots$ denote terms arising from bringing at least the
coefficient in $aT$ to the left of $X^r,$ incurring additional
derivatives of course on the coefficient $a.$ However all of this
is linear. The non-linear phenomena occur in the last term, where
$a_u = a(x,u)$ is differentiated. But this proceeds as
before (cf. (\ref{[Tr,P]})): letting, for instance, 
$b^{(s)}(x,u)$ denote derivatives of the function $b$ in {\it
its} arguments, all derivatives of the solution $u$ being split
off to the right,
\begin{equation}{[a_u, X^r]w \over r!}
=\sum_{{{r\geq r'-r''\geq 1}\atop \sum_{j=1}^{r'-r''} r''_j =
r''}\atop {(\sum_{j=1}^{r'-r''} (r''_j+1) =
r'})}{a_u^{(r'-r'')}\over
{(r'-r'')}!}\;{X^{r''_1}{u'}\over r''_1!}\cdots
{X^{r''_{r'-r''}}{u'}\over r''_{r'-r''}!}\;{X^{r-r'}w\over
(r-r')!}.
\end{equation}
So, all together,
(\ref{exp:[P, Psi X^r]-1}) becomes:
\begin{multline}
\label{exp:[P, Psi X^r]-2}[ P,\Psi X^r]v \sim 
2\,a_u\,
\Psi 'XX^rv + 
a_u\, \Psi ''X^rv + 
\,\underline{r}\;a_u\,\Psi \,a\,TX^rv+\cdots 
\\
+\Psi\,r!\sum_{{{r\geq r'-r''\geq 1}\atop \sum_{j=1}^{r'-r''}
r''_j = r''}\atop {(\sum_{j=1}^{r'-r''} (r''_j+1) =
r'})}{a_u^{(r'-r'')}\over
{(r'-r'')}!}\;{X^{r''_1}{u'}\over r''_1!}\cdots
{X^{r''_{r'-r''}}{u'}\over r''_{r'-r''}!}\;{X^{r-r'}X^2v\over
(r-r')!}.
\end{multline}

Now once we specialize $v$ to $u,$ we will take the $H^s$
norm of everything and use the property that this space is an
algebra for our choice of $s.$  The function $\Psi$ in the
product on the right has served its purpose, and we will
eventually introduce a new localizing function for each term in
the product (except the coefficient, which will just be estimated),
though at most one of these terms can have `order' even half of
$r$ and the rest will be handled inductively.  

\section{Local Regularity in High Powers of $X;$ new
localizing functions}
\renewcommand{\theequation}{\thesection.\arabic{equation}}
\setcounter{equation}{0}
\setcounter{theorem}{0}
\setcounter{proposition}{0}  
\setcounter{lemma}{0}
\setcounter{corollary}{0} 
\setcounter{definition}{0}
\setcounter{remark}{1}

More precisely, we restate (\ref{exp:[P, Psi X^r]-1}) after specialization
and introduction of the $H^s$ norm:
$${\|[ P,\Psi X^r]u\|_s\over
r!}
\leq  C_u\{
{\|\Psi 'X^{r+1}u\|_s\over r!} + 
{\|\Psi ''X^ru\|_s\over r!} + 
\underline{r}{\|\Psi TX^ru\|_s\over r!} +\cdots \}
$$
$$+\sum_{{{r\geq r'-r''\geq 1}\atop \sum_{j=1}^{r'-r''}
r''_j = r''}\atop {(\sum_{j=1}^{r'-r''} (r''_j+1) =
r'})}\|\Psi {a_u^{(r'-r'')}\over
{(r'-r'')}!}{X^{r''_1}{u'}\over r''_1!}\cdots
{X^{r''_{r'-r''}}{u'}\over r''_{r'-r''}!}\;{X^{r-r'}X^2u\over
(r-r')!}\|_{H^s}.
$$
We treat the functions $X^2u$ and $u'$ similarly
- they are equivalently handled by the {\it a priori} estimate - and
{\it for convenience only} we suppose that the term with $X^2u$
is of highest order - i.e., $r-r'\geq r_j''\;\;\forall j.$ Noting
that ${\hbox{ supp}}\Psi \Subset \mathcal{U}_{1/r}$
and bounding the norm of the coefficients by $C^{r'-r''},$
%split off last factor
\begin{multline}
\label{est:||[P, Psi X^r]||_s} 
{\|[ P,\Psi X^r]u\|_s\over
r!}
\leq  C_u\{
{\|\Psi 'X^{r+1}u\|_s\over r!} + 
{\|\Psi ''X^ru\|_s\over r!} + 
\underline{r}{\|\Psi TX^ru\|_s\over r!} +\cdots \}
\\
+\sum_{{{r\geq r'-r''\geq 1}\atop \sum_{j=1}^{r'-r''}
r''_j = r''}\atop {(\sum_{j=1}^{r'-r''} (r''_j+1) =
r'})}C^{r'-r''+2}\left \|
{\prod}_{j=1}^{r'-r''}{X^{r''_j}{u'}\over r''_j!}\right
\|_{H^s(\mathcal{U}_{1/r})}\,
\left\|{\Psi X^{r-r'}X^2u\over
(r-r')!}\right\|_{H^s}.
\end{multline}
Again, we note that the number of terms in the product is $r'-r''$
and hence the constant arising from the algebraicity of
$H^s$ will be absorbed with the analyticity constant for
the coefficients $a_u.$  Thus we restate (\ref{est:||[P, Psi
X^r]||_s}) with this observation, writing $\Psi X^2$ in place of
$X^2\Psi$ on the left, modulo terms on the right, and associating
powers of $r$ with derivatives of $\Psi$ or with powers of $T,$
and taking $Pu=0:$%
\begin{multline}
\label{onefullpass} 
{\|  \Psi X^2 X^r u\|_s\over r!}
\leq C_u\{
\sum_{j=1}^2{{1\over r^j}\|\Psi^{(j)}X^{r+2-j}u\|_s\over (r-j)!} + 
{{1\over r}\|\Psi TX^ru\|_s\over (r-2)!} +\cdots \}
\\
+\sum_{{{r\geq r'-r''\geq 1}\atop \sum_{j=1}^{r'-r''}
r''_j = r''}\atop {(\sum_{j=1}^{r'-r''} (r''_j+1) =
r'})}C_h\left \|
{\prod}_{j=1}^{r'-r''}{X^{r''_j}{u'}\over
r''_j!}\right
\|_{H^s(K)}\,
\left\|{\Psi X^2X^{r-r'}u\over
(r-r')!}\right\|_{H^s}.
\end{multline}

As we iterate the terms on the right without $T,$
the order will drop and we will control the
coefficients and the sum below. The term with $T$ is
slightly different, but we may always write,  
$${{1\over r}\|\Psi TX^ru\|_s\over (r-2)!} = 
{{1\over r}\|\Psi X^2 X^{r-2} Tu\|_s\over (r-2)!}.
$$ 
and reapply (\ref{onefullpass}) with $Tu$ in place of $u$ 
but with $r$ decreased by two. Thus we 
gradually increase the number of
$T$ vector fields, with $T^\sigma$ being
balanced by $1\over \sigma!!$ before the norm, where
$$\sigma!! = \sigma(\sigma-2)(\sigma-4)\ldots,$$
preserving the balance between remaining powers of $X$
and the large factorial in the denominator, and using up two
$X$'s for each
$T$ until there are essentially only powers of
$T,$ a situation we have treated above. (Of course there
will be a `zig-zag' effect - sometimes pairs of $X$'s
will generate a $T$ and other times the $X'$s will
differentiate the coefficients and produce the terms at
the end of (\ref{onefullpass}) above, so both effects will
be combined.)  

And as with the estimates of pure
$T$ derivatives above, iterating the `principal' term
(here the last one - the one with
$X^2 X^{r-r'} u$) will lead to a sum with the same
bounds for the number of its terms (cf.
(\ref{partitions})), and with one new norm of derivatives
of a coefficient
$a_u.$ Even when $\Psi$ has not been differentiated, it will be
prudent to change to a new localizing function, one better geared
to the number of derivatives appearing under the 
norm. For
there are fewer derivatives now, and it would create
significant difficulties to have
$\Psi'$ contribute a factor of $r$ when the denominator
contains only $(r-r')!$ for
rather general
$r'.$

\section{The Localizing Functions}
\renewcommand{\theequation}{\thesection.\arabic{equation}}
\setcounter{equation}{0}
\setcounter{theorem}{0}
\setcounter{proposition}{0}  
\setcounter{lemma}{0}
\setcounter{corollary}{0} 
\setcounter{definition}{0}
\setcounter{remark}{1}

The first localizing function, $\Psi = \Psi_r,$
satisfies:
\begin{equation}
\label{def:Psi_r} 
\Psi_r\equiv 1  \mbox{ on }\ {\mathcal U}_0\ ,\
\Psi_r\in
C_0^{\infty}({\mathcal U}_{1/r})\ ,
\ |\Psi_r^{(k)}|\leq c^kr^{k}, \;k\leq p(s),\end{equation}
(cf. (\ref{Psirhosigma2})), where we have set, for $a\geq 0:$
\begin{equation}
\label{def:Ua} 
{\mathcal U}_a=\{(x,t)\in{\mathcal U}_1\ :\
\mbox{dist}((x,t),{\mathcal U}_0)<a(\text{dist}(\mathcal U_0,
\mathcal U_1^c)\}\ .\end{equation} 

When the first  
localizing function needs to be replaced but, say, 
 $\tilde{r}$ derivatives of $u$ remain to be estimated, 
we shall localize it with a function identically
equal to one on $\mathcal U_{1/r},$ the support of $\Psi_r$
but dropping to zero in a band of width
${1\over \tilde{r}}\times (1-{1\over
r})(\text{dist}(\mathcal U_0,
\mathcal U_1^c)) = {1\over \tilde{r}}$ times the
remaining distance
to the complement of 
$\mathcal U_1,$ i.e., supported in 
\begin{equation}
\label{supports1}
\mathcal U_{{1\over r}+({1\over
\tilde{r}})(1-{1\over r})}=\mathcal U_{{1\over r}+{r-1\over
r\tilde{r}}}=\mathcal U_{1-(1-{1\over r})(1-{1\over
\tilde{r}})}.
\end{equation}
We shall denote such a function by $_{1\over r}\Psi
_{\tilde{r}}$ 
That is,
${}_\rho\Psi_\sigma$ satisfies:
\begin{equation}
\label{Psirhosigma1} 
{}_\rho\Psi_\sigma \equiv 1 \mbox{ on }
{\mathcal U}_{\rho}, \qquad
{}_\rho\Psi_\sigma \in
C_0^{\infty}({\mathcal{U}}_{\rho + {1\over
\sigma}(1-\rho)}\Subset {\mathcal
U}_1).
\end{equation}  
Derivatives of ${}_\rho\Psi_\sigma$ satisfy, with
universal constant $C$: 
\begin{equation}
\label{Psirhosigma2}
|D^k \left({}_\rho\Psi_\sigma\right)| \leq C^k\left({\sigma \over
1-\rho}\right)^k,
\;k\leq p(s).
\end{equation}
uniformly in $\rho, \sigma,$ where $p(s)$ will be a
small number depending on the $s$ necessary to make
$H^s$ an algebra in the given dimension. Of course
any other (fixed) bound for
$k$ would do.  
\section{Taking a localizing function out of the norm}
\renewcommand{\theequation}{\thesection.\arabic{equation}}
\setcounter{equation}{0}
\setcounter{theorem}{0}
\setcounter{proposition}{0}  
\setcounter{lemma}{0}
\setcounter{corollary}{0} 
\setcounter{definition}{0}
\setcounter{remark}{1}

While it is true that we could just write
$\|\Psi w\|_s\leq c\|\Psi\|_s\|w\|_s,$ for
$s>1,$ to do so would incur at least two derivatives on
$\Psi$ with no gain on $w.$ To avoid this
difficulty, we use the following finer estimates of the $H^s$ norm
of product of functions.
\begin{proposition}\label{{stos+infinity0}}
   If $\Psi,\tilde\Psi$ are two smooth, compactly
supported functions with $\tilde\Psi\equiv 1$
on supp $\Psi$ then for every $s\leq p \in \Z^+,$
\begin{equation} \label{{stos+infinity1}}\|\Psi D^pu\|_s\leq 
C_{s, supp\,\Psi}^2\sup_{q\leq
s}\|D^q\Psi\|_{L^{\infty}}\|\tilde\Psi
D^{p-q}u\|_s 
\end{equation}
and
\begin{equation} \label{stos+infinity2} \|\Psi D^pu\|_s\leq
C_{s, supp\,\Psi}^2\sup_{q\leq
s}\|D^q\Psi\|_{L^{\infty}}\|D^{p-q}u\|_{H^s(supp\
\Psi)}\ .\end{equation}
\end{proposition}

Thus removing a localizing function from an $H^s$ norms,
while incurring up to $2$ derivatives on it,
does not increase the total number of derivatives being
measured, and thus should have minimal impact on the
estimates.

Next, we need to confront the effect of these few derivatives on
a localizing function $\Psi,$ which may have been chosen with a
high number of derivatives ($r$ of them) on $u$ in
mind, and hence which adds a factor of $r$ each time a derivative
lands on it, when the factorial in the denominator of the
corresponding term may be far smaller, e.g., $(r-r')!$ for
relatively large $r'.$ 

There are in fact several ways to handle
this; one is to emphasize that at the level of
(\ref{bracketwithD^r}) one could endeavor to keep two
derivatives to the left of the big bracket whenever possible so
that using Proposition \ref{{stos+infinity0}}, those derivatives
would serve to bring the Sobolev norms on the right up to $H^s$
again, or one can proceed as follows, the method used in the
second author's earlier work
\cite{Tartakoff-Zanghirati2003}: since the number of terms in the
product in (\ref{est:||[P, Psi X^r]||_s}) is $r'-r'',$ with
\begin{equation} 
r=(r-r')+\sum_{j=1}^{r'-r''} (r''_j+1),
\hbox{ with } r-r' \geq \,\max \,\{r_j''+1\}
\end{equation}
and so
\begin{equation}
(r-r')(r'-r''+1)\geq r
\end{equation}
or
\begin{equation}\label{r,r-r'}
\left({r\over r-r'}\right)^k \leq (r'-r''+1)^k, \,\forall k
\end{equation}
a relationship we will use only for small values of $k$ but note
that this factor, $(r'-r''+1)^k,$ can be absorbed in the bound of
derivatives of the coefficients $a_u$ in (\ref{est:||[P, Psi X^r]||_s}).

The first time we remove a localizing function from an $H^s$
norm, in (\ref{onefullpass}) for instance, the couple of
derivatives that will fall on $\Psi$ will produce powers of $r$ in
view of (\ref{Psirhosigma2}), since initially $\rho=0.$ These
will be balanced against $(r-1)!$ thanks to (\ref{r,r-r'}) with 
small powers of $(r'-r''+1),$ increasing $C_h$ slightly in
(\ref{onefullpass}). We will see at the very end that the slightly
different denominators in (\ref{Psirhosigma2}) will make little
difference in the bounds. 

Furthermore, upon the next iteration of  (\ref{onefullpass}), the new right
hand side {\it will have the same form}. That is, there will again
be a product of lower order terms (the same ones plus new ones),
a second copy of $a_u$ with derivatives which will give possibly
another constant,  
$C$ in front of the supremum and another copy of $C_h$ to its appropriate
power, though these constants pass into the norms of the
corresponding terms, just as in the treatment of powers of $T$
above. But notice that the number of terms in the
product increases at each pass (to at most
$r$) and the the order of the top order term decreases. Thus this
sequence of constants will not contribute in the end more than
$C^r,$ which is also to be expected. 

Handling the sum is as before as well, and we will not
comment on it further except to recall (\ref{partitions}). 

When
the last term on the right no longer has maximal order, we turn our
attention to any of the other terms of highest order and proceed as
before. The factorials have been adjusted so that the behavior that
will in the end guarantee analyticity is that 
$${|| \Psi X^2 X^ru||_{H^s}\over (r-1)!}\leq
C^{r+1} + C^{r/2}{||X^2 T^{r/2}
u||_{H^s(\mathcal{U}_1)}\over (r/2)!}$$  
which will be the
evident outcome of the repeated iterations of (\ref{onefullpass})
taking the precise localizing functions into account and which,
together with the previous results on (nearly) pure $T$
derivatives will complete the proof. 

We should remark at the end that what was true for the first localizing
function, namely (\ref{r,r-r'}), will be a little different on
the next pass, since the next localizing function may bring not a
factor of
$r-r'$ with each derivative it receives but rather the factor (cf.
(\ref{Psirhosigma2})) 
$${r-r'\over 1-{1\over r}}= (r-r')\left({r\over r-1}\right)$$
so that, passing from $r-r'$ to $r-r'-t'$ we encounter instead of
just 
$${r \over {r-r'}}\leq r'-r''+1$$
an extra factor of $r/r-1,$ possibly to the $p(s)$-th power; and this
may  keep occurring as the order of the leading term keeps
decreasing. For instance, after a few iterations, the analogous
`extra' factors from (\ref{Psirhosigma2})  will be
$$\left({r\over r-1}\right)\left({r-r_1\over
r-r_1-1}\right)\left({r-r_1-r_2\over
r-r_1-r_2-1}\right)\ldots
$$  or even the $p(s)$-th power of such a product. But there cannot
be more than
$r$ terms in the product and each factor is far less than $2,$
leading to an easily acceptable constant $C^r$ in the end. 

The same procedure works at any stage. We have already seen
that expanding the term of highest order leads to a new product,
but of the {\it same form} with one new norm of derivatives of a
coefficient
$a_u,$ and the total number of terms, as with the $T$ derivatives,
never exceeds
$4^r,$ which is certainly acceptable; and the factor $(r'-r''+1)^k$
just above is immediately attached to the $a_u^{(r'-r'')}$ which
occurs with that product (cf.(\ref{exp:[P, Psi X^r]-2})). 

This means that we may remove localizing functions from the
$H^s$ norms easily and replace them with new localizing
functions, identically one on the support of the old one and
supported in a larger open set such that a derivative of the new
function is proportional to the number of derivatives still to be
estimated in that term in the sense of
(\ref{Psirhosigma2}). And while there will appear a number of
copies of the (analytic) coefficients $a,$ the sum of the number
of derivatives they receive, and the powers of the corresponding
constants arising from the algebraicity of $H^s,$ is equal to the
total decrease in derivatives on the terms of highest order taken
step by step, which is also reflected in the number of norms in
the product. Thus the total number of derivatives appearing on
coefficients will be, in the end, equal to the total number of
terms in the product of norms - and since each contains a copy of
$u$ with one or two derivatives, this number is comparable to
$r.$ Thus this product will be bounded by $C_{u,h}^r$ for
suitable $C_{u,h}$ depending only on the first couple of
derivatives of $u$ and on the coefficients $a_u$ (and the
dimension and the initial open sets). 

We are not quite home. For high $X$ derivatives, in addition to
being `used up' as above, will also flow to half as many
$T$ derivatives, though in
$\mathcal{U}_1,$ due to the bracketting $[X,X]=T$ (cf.
(\ref{[X,X]=aT})), and the number of terms and the sums proceed
exactly as in the estimation of $T$ derivatives above, in ways
that have nothing to do with the local versus global behavior.
There appear new norms of derivatives of the coefficient
functions, exactly as before, and one slightly new feature which is
the mixture between $X$ and
$T$ derivatives which is inevitable but has been seen before in
many of the authors' earlier works.

\end{document}